\documentclass[12pt]{article}
\usepackage{amssymb, amsmath,amsthm,amsfonts,a4,bezier}
\usepackage{hyperref}
\usepackage{amscd}
\usepackage{graphpap}
\usepackage[pdftex]{color,graphicx}
\usepackage{color,graphicx}
\usepackage{mathpazo}
\usepackage{setspace}
\usepackage{comment}
\usepackage{mathtools}
\usepackage{relsize}

\DeclareGraphicsExtensions{.jpg,.pdf,.eps,.png}
\newtheorem{theo}{Theorem}

\newtheorem{cor}{Corollary}

\newtheorem{Remark}{Remark}
\newtheorem{Conjecture}{Conjecture}
\renewcommand{\o}{{\omega}}
\newcommand{\z}{{\mathbf {z}}}

\renewcommand{\t}{{\mathbf {t}}}
\newcommand{\W}{{\mathbf {W}}}
\renewcommand{\u}{{\mathbf {u}}}
\newcommand{\Z}{{\mathbb {Z}}}

\newcommand{\E}{{\mathbb {E}}}
\newcommand{\F}{{\mathcal {F}}}

\begin{document}
\thispagestyle{empty}
\baselineskip=28pt
\vskip 5mm
\begin{center} {\LARGE{\bf Dimensional Crossover in Anisotropic Percolation on $\mathbb{Z}^{d+s}$}}
\end{center}

\baselineskip=12pt
\vskip 10mm

\begin{center}\large
R\'{e}my Sanchis\footnote{Departamento de Matem\'{a}tica, Universidade Federal de Minas Gerais, Brazil, \href{mailto:rsanchis@mat.ufmg.br}{rsanchis@mat.ufmg.br}} 
and Roger W. C. Silva 
 \footnote{Departamento de Estat\'{\i}stica,
 Universidade Federal de Minas Gerais, Brazil, \href{mailto:rogerwcs@ufmg.br} {rogerwcs@ufmg.br}}
\\ 


\end{center}

\begin{abstract} We consider bond percolation on $\Z^d\times \Z^s$ where edges of $\Z^d$ are open with probability $p<p_c(\Z^d)$ and edges of $\Z^s$ are open with probability $q$, independently of all others. We obtain  bounds for the critical curve  in $(p, q)$, with $p$ close to the critical threshold $p_c(\Z^d)$. The results are related to the so-called dimensional crossover from $\Z^d$ to $\Z^{d+s}$.\\


\end{abstract}

\onehalfspacing
\section{Introduction and Results}
\subsection{Background}\label{background}

\hspace{0.5cm}In this note we consider an anisotropic bond percolation process on the graph \linebreak
$(\Z^{d+s},\mathcal{E}(\Z^{d+s}))$, where $\mathcal{E}(\Z^{d+s})$ is the set of edges between nearest neighbors of $\Z^{d+s}$. With some abuse of notation, we call this graph  $\Z^{d+s}=\Z^d\times \Z^s$. An edge of $\Z^{d+s}$ is called a $\Z^d$- edge (respectively a $\Z^s$- edge) if it joins two vertices which differ only in their $\Z^d$ (respectively $\Z^s$) component. Given two parameters $p,q\in[0,1]$, we declare each $\Z^d$- edge open with probability $p$ and each $\Z^s$- edge open with probability $q$, independently of all others. This model is described by the probability space $(\Omega, \F, \mathbb{P}_{p,q})$ where
$\Omega = \{0,1\}^\mathcal{E}$, $\F$ is the $\sigma$-algebra generated by the cylinder sets in $\Omega$ and $\mathbb{P}_{p,q} =\prod_{e\in\E}\mu(e)$, where $\mu(e)$ is Bernoulli measure with parameter $p$ or $q$ according to $e$ been a $\Z^d$- edge or a $\Z^s$- edge respectively. 

Given two vertices $u$ and $v$, we say that $u$ and $v$ are connected in the configuration $\o$  if there exists an open path starting in $u$ and ending in $v$. We denote the event where $v$ and $u$ are connected by $\{\o\in\Omega:v\leftrightarrow u\mbox{ in $\o$}\}$ and write $C(\o)=\{u\in \Z^d\times\Z^s:u\leftrightarrow 0\mbox{ in $\o$}\}$ for the open cluster containing the origin. We write $\theta(p,q)=\mathbb{P}_{p,q}(\o\in\Omega:|C(\o)|=\infty)$ for the main macroscopic function in percolation theory and denote the mean size of the open cluster by $\chi(p,q)=\E_{p,q}(|C(\o))|)$. Whenever necessary, we will write  $\chi_d(p)$ and $p_c(d)$  in reference to the expected cluster size and critical threshold on $\Z^d$ with a single parameter $p\in(0,1)$. For a thorough background in independent percolation, we refer the reader to \cite{GRI}.

By a standard coupling argument we have that $\theta(p,q)$ is a monotone non decreasing function of the parameters $p$ and $q$. This allows us to define the function $q_c:[0,1]\rightarrow [0,1]$, where
\begin{equation}\label{qcrit}q_c(p)=\sup\{q:\theta(p,q)=0\}.
\end{equation}
We observe that $q_c(0)=p_c(s)$ and $q_c(p)=0$ for $p\geq p_c(d)$. Moreover, $q_c(p)$ is a non increasing function of the parameter $p$. Our main interest is the behavior of $q_c(p)$ for values of $p$ less or equal than $p_c(d)$, in particular when $p\uparrow p_c(d)$.

\subsection{Dimensional Crossover}
\hspace{0.5cm}Consider a bond percolation process on the graph $\Z^d\times\Z^s$, as defined in Section \ref{background}. The problem we address in this note is the following: what is the behavior of the critical curve $q_c(p)$ as $p\uparrow p_c(d)$? 

A similar model was  investigated in \cite{CLS}, where the authors study
anisotropic percolation on the slab $\Z^2\times\{0,\dots,k\}$. In their model, vertical and horizontal edges are open with probability $p$ and $q$ respectively, independently of all others. They obtain bounds for the critical curves for these models and establish their continuity and strict monotonicity. Similar anisotropic ferromagnectic models have also been considered in the mathematical literature, see the works \cite{FMMPV}, \cite{FMMPV2} and  \cite{MPS}. 

A central open problem in percolation theory is the existence and determination of critical exponents. For instance,  consider the isotropic percolation process with parameter $p$ on $\Z^d$. Quantities such as $\chi_d(p)$ are believed to diverge as $p\uparrow p_c(d)$ in the manner of a power law in $|p-p_c(d)|$, whose exponent is called a critical exponent (see Chapter 9 in \cite{GRI} for details). That is, it is believed that there exists a $\gamma(d)>0$ such that 
\begin{equation}\label{crit}
\chi_d(p)\approx |p-p_c(d)|^{-\gamma},
\end{equation}
when $p\uparrow p_c(d)$. Here the relation $a(p) \approx b(p)$ means {\textit log} equivalence, i.e., $\frac{\log a(p)}{\log b(p)}\rightarrow 1$ when ${p\uparrow p_c(d)}$.

We introduce another critical exponent in the following way: consider the function defined in Equation (\ref{qcrit}). Is it true that there exists a constant $\psi>0$ (which depends on $d$) such that $q_c(p)\approx |p-p_c(d)|^{\psi}$ when $p\uparrow p_c(d)$? This problem arises in the physics literature as the dimensional crossover problem, and $\psi$ is called the crossover exponent. The term crossover is related to the study of percolative systems on  $(d+s)$-dimensional lattices, where the $d$-dimensional parameter $p$ is close to $p_c(d)$ from below and the $s$-dimensional parameter $q$ is small. 

It is commonly suggested in the physics literature that $\psi(d)$ and $\gamma(d)$ exist and are equal. Let us highlight a few papers that discuss this matter.  In \cite{GCGR} the authors study anisotropic bond percolation on $\Z^3=\Z^2 \times \Z$ . Here $\Z^2$- edges are open with probability $p$ and $\Z$- edges are open with probability $q=Rp$, where $R$ is the anisotropy parameter. By means of a simulation, the authors estimate $\psi(2)$ by $2.3 \pm 0.1$, which is compatible with the critical exponent $\gamma(2)$, whose estimate given in \cite{S} is $\frac{43}{18}$.
In \cite{RS}, Redner and Stanley consider a percolation model on $\Z^d=\Z^{d-1}\times\Z$ where $\Zö^{d-1}$- edges are open with probability $p$ and $\Z$- edges parallel to $z$ are open with probability $q= Rp$.  Through simulated results the authors show that in the limit  $1/R\rightarrow 0$,  the crossover exponent $\psi$ is equal to $1$ for all $d$. In the opposite limit $R\rightarrow 0$, their analysis suggests that $\psi(d-1)\neq\gamma(d-1)$. This result was later contradicted by Redner and Coniglio \cite{RC}, where the authors argue the opposite relation, that is,  $\psi(d-1)=\gamma(d-1)$. The interested reader can consult the works \cite{BMT},\cite{DB},\cite{FS} and \cite{M} for more results in this direction.

 Based on the results mentioned above, we state the following:

\begin{Conjecture}\label{crit2} There exists a critical exponent $\psi=\psi(d)>0$ such that
$$q_c(p)\approx |p-p_c(d)|^{\psi}.$$ Besides, if $\gamma(d)$ exists, then $\psi(d)=\gamma(d)$.
\end{Conjecture}

In Section \ref{result} we give a complete answer to this conjecture in the case $d=1$ and a partial answer for general $d$. We remark that the crossover critical exponent, if it exists, depends only on the dimension $d$ and not on $s$.

\subsection{Results} \label{result}

\hspace{0.5cm} We start with the following theorem, which gives a lower bound for the critical curve $q_c(p)$.

\begin{theo}\label{lemma1} Consider a bond percolation process on $\Z^d\times\Z^s$ with parameters $(p,q)$, $p<p_c(d)$. If the pair $(p,q)$ satisfies $$q<\frac{1}{2s\chi_d(p)}\,,$$ then there is a.s. no infinite open cluster in $\Z^{d+s}$.
\end{theo}

Unfortunately, we were not able to obtain an upper bound for $q_c(p)$ that would work for any $d$ and $s$. Nevertheless, in the particular case $d=1$, having full control of the one dimensional cluster  enables us to build up the open cluster of the origin in  $\Z^{1+s}$ by tuning the parameters appropriately.   We adopt notation from Section \ref{background} and say $\Z$- edges are open with probability $0<p<1$, independently of all other edges and $\Z^s$- edges are open with probability $0<q<1$, also independently of all other edges. The next theorem gives an upper bound for the critical curve $q_c(p)$. 

\begin{theo}\label{lemma2} Consider a bond percolation process on $\Z\times\Z^{s}$, $s>1$, with parameters $(p,q)$.  
  There exists an $\alpha>0$ such that if  $p$ is sufficiently close to $1$ and $q>\alpha \frac{1-p}{1+p}$, then there is a.s. an infinite open cluster in $\Z^{d+1}$.
\end{theo} 
Since $\chi_1(p)=\frac{1+p}{1-p}\,,$ an immediate consequence of Theorems \ref{lemma1} and \ref{lemma2} is    
 
\begin{cor}\label{theo1}
Consider a bond percolation process on $\Z\times\Z^{s}$, $s>1$, with parameters $(p,q)$. Then $$\frac{1}{2s\chi_1(p)}\leq q_c(p)\leq \frac{\alpha}{\chi_1(p)}\,,$$ for some $\alpha>0$ and $p$ sufficiently close to $1$.
\end{cor}

\begin{Remark}  By Theorem \ref{lemma1}, the lower bound for $q_c(p)$ in Corollary \ref{theo1} remains valid for any $p\in[0,1)$.
\end{Remark}
\begin{Remark} The case $s=1$ has been fully understood in \cite{K}. It was shown that the critical curve is of the form $\psi(p,q)=p+q-1$, that is, if $\psi(p,q)<0$, then there is no infinite open cluster a.s.,  whereas if $\psi(p,q)>0$, then there exists an infinite open cluster a.s..
\end{Remark}

Another straightforward consequence of Theorems \ref{theo1} and  \ref{lemma2} is 
\begin{cor}\label{corexp} We have the following relations for the critical exponents $\gamma$ and $\psi$:

\begin{enumerate}
\item $\psi(1)$ exists and is equal to $\gamma(1)=1$. Hence Conjecture \ref{crit2} is true in the case $d=1$.
\item If $\gamma(d)$ and $\psi(d)$ exist, then $\psi(d)\leq \gamma(d)$ for all $d$.
\end{enumerate}
\end{cor}

We end this section with the following conjecture, which together with Theorem \ref{lemma1} and the validity of Equation (\ref{crit}) implies that Conjecture \ref{crit2} is true in general. In particular it establishes that, if the critical exponents $\psi (d)$ and $\gamma(d)$ exist, then they are equal.

\begin{Conjecture}\label{conj2} Consider a bond percolation model on $\Z^d\times\Z^{s}$,  with parameters $(p,q)$, $p<p_c(d)$. Then $$q_c(p)\leq\frac{\beta}{\chi_d(p)}\,,$$ for some $\beta>0$.
\end{Conjecture}

\section{Proofs }
\subsection{Proof of Theorem \ref{lemma1}}
\hspace{0.5cm} Let us first introduce some notation. For any $(\u,\t)\in\Z^{d+s}$, $\u$ and $\t$ will always denote the $\Z^d$ and $\Z^s$ components of $(\u,\t)$ respectively. We shall also consider sequences of points $\{\u_n,\, n\in\Z^+\}$ and $\{\t_n,\, n\in\Z^+\}$ in $\Z^d$ and $\Z^s$ respectively. We let $\mathbf{0}_d$  denote a $d$-dimensional vector of zeros. Finally, given $\u\in\Z^d$ and $\mathbf{s},\t\in \Z^s$ such that $\delta(\mathbf{s},\t)=1$, write $e_{(\u,\t),(\u,\mathbf{s})}\in \mathcal{E}(\Z^s)$  for the edge with endpoints $(\t,\u)$ and $(\u,\mathbf{s})$. Here, $\delta(\mathbf{x},\mathbf{y})=\sum_{i=1}^s(x_i-y_i)$.

The following argument is based on techniques developed in \cite{GN}. We will bound $\chi(p,q)$, the expected size of the infinite open cluster in $\Z^{d+s}$. Clearly,
\begin{equation*}\chi(p,q)=\sum_{(\u,\t)\in \Z^{d+s}}\mathbb{P}_{p,q}(\mathbf{0}_{d+s}\leftrightarrow (\u,\t)).
\end{equation*} 
Now, the event $\{\mathbf{0}_{d+s}\leftrightarrow (\u,\t)\}$ occurs if and only if there exist 
sequences $\hat{\u}_n=(\u_0,\dots\u_n)$ and $\hat{\t}_n=(\t_0,\dots\t_n)$ such that $\u_j\in\Z^d$, for all $0\leq j\leq n$, $\t_i\in\Z^s$ for all $1\leq i\leq n$ and $\delta(\t_{i},\t_{i-1})=1$, with the following property: we start at $(\mathbf{0}_d, \t_0)=(\mathbf{0}_d, \mathbf{0}_s)$ and then we connect      
$(\mathbf{0}_d, \t_0)$ to $(\u_0,\t_0)$ using only $\Z^d\times \{\t_0\}$- edges. Next, we move from $(\u_0,\t_0)$ to $(\u_0,\t_1)$ by a connection that uses a single open $\Z^s$- edge. We proceed looking for a connection from $(\u_0, \t_1)$ to $(\u_1,\t_1)$ that uses only $\Z^d\times\{\t_1\}$- edges that are not used to connect $(\mathbf{0}_d, \t_0)$ to $(\u_0,\t_0)$. At the $m-th$ step we move from $(\u_{m-1},\t_{m-1})$ to $(\u_{m-1},\t_{m})$ by a connection that uses a single open $\Z^s$- edge which is not used in any of the previous steps of the construction. Next, we connect $(\u_{m-1},\t_m)$ to $(\u_{m},\t_m)$ using only $\Z^d\times\{\t_m\}$- edges, with the  constraint that no $\Z^d$-  edges used to connect $(\u_{l-1}, \t_l)$ to $(\u_l, \t_l)$, $1\leq l\leq m-1$, and $(\mathbf{0}_d, \t_0)$ to $(\u_0,\t_0)$ are used to connect $(\u_{m-1},\t_m)$ to $(\u_m,\t_m)$.  We continue this procedure until it reaches $(\u_n,\t_n)=(\u,\t)$.   

Given $\hat{\t}_n=(\t_0,\dots,\t_n)$ and $\hat{\u}_n=(\u_0,\dots,\u_n)$, consider the sequence of increasing events $\{A_i\}_{i=0}^n$, where, for  $0\leq i \leq n-1$,
$$A_i=\{\{(\u_{i-1},\t_{i})\leftrightarrow (\u_{i},\t_{i})\mbox{ on $\Z^d\times\{\t_i\}$}\}\cap \{e_{(\u_i,\t_i),(\u_i,\t_{i+1})}\mbox{ is open}\}\},$$
and $$A_n=\{\{(\u_{n-1},\t_{n})\leftrightarrow (\u_{n},\t_{n})\mbox{ on $\Z^d\times\{\t_m\}$}\}\},$$ with the convention that $\u_{-1}=\mathbf{0}_d$.
As a consequence of the above construction we have
$$\{\mathbf{0}_{d+s}\leftrightarrow (\u,\t)\}=\mathlarger{\bigcup_{n\geq 0}\bigcup_{\substack{\hat{\u}_n:\u_n=\u\\ \hat{\t}_n:t_n=\t}}\{A_0\circ A_1\circ \dots \circ A_n\}},$$
where, $A\circ B$ denotes disjoint occurrence of events $A$ and $B$. By the BK inequality (see \cite{BK}) we have
$$\chi(p,q)\leq \mathlarger{\sum_{(\u,\t)\in \Z^{d+s}}\sum_{n\geq 0}\sum_{\substack{\hat{\u}_n:\u_n=\u\\ \hat{\t}_n:t_n=\t}}\prod_{i=0}^n\mathbb{P}_{p,q}(A_i)},$$ or equivalently
$$\chi(p,q)\leq \mathlarger{\sum_{n\geq 0}\sum_{\hat{\t}_n}\sum_{\hat{\u}_n}\prod_{i=0}^n\mathbb{P}_{p,q}(A_i)},$$ where the
final two summations are over all appropriate sequences $\hat{\u}_n=(\u_0,\dots,\u_n)$ and $\hat{\t}_n=(\t_0,\dots,\t_n)$, where $(\u_0,\t_0)=\mathbf{0}_{d+s}$. Clearly,
$$\sum_{\u_n\in\Z^d}\mathbb{P}_{p,q}(A_n)\leq \chi_d(p),$$ and, since the events $\{(\u_{i-1},\t_{i})\leftrightarrow (\u_{i},\t_{i})\mbox{ on $\Z^d$}\}$ and $\{e_{(\u_i,\t_i),(\u_i,\t_{i+1})}\mbox{ is open}\}$ are independent, we obtain
$$\sum_{\substack{\u_k\in\Z^d\\ \t_k\in\Z^s}}\mathbb{P}_{p,q}(A_k)\leq 2sq\chi_d(p),\,\,\,\,\,\,\,\forall \, k=0,\dots,n-1.$$ Proceeding inductively we obtain

$$\chi(p,q)\leq\sum_{n=0}^{\infty}(2s)^n\chi_d(p)^{n+1}q^n.$$
If $p<p_c(d)$, then $\chi_d(p)<\infty$ a.s. and the result follows.


\subsection{Proof of Theorem \ref{lemma2}}
\hspace{0.5cm}We shall construct an independent site percolation process in $\Z^{s}$ which is induced by the bond percolation process in $\Z^{1+s}$. We show that, under the appropriate hypothesis, site percolation occurs in $\Z^{s}$, and therefore, by stochastic dominance, bond percolation will occur in $\Z^{1+s}$.

Let $\u=(u,u_1,\dots,u_s),\z=(z,z_1,\dots,z_s)\in\Z^{1+s}$. Given a configuration $\omega\in\Omega$, we say that $\u$ is updownwards connected 
  to $\z$ in the configuration $\omega$ if $u_i=z_i$ for all $i=1,\dots,s$, and there exists a open path using only $\Z\times\{u_1,\dots,u_s\}$- edges starting at $\u$ and ending at $\z$.
  Denote this event by $\{\omega\in\Omega: \u\updownarrow 
  \z\mbox{ in }\omega\}$. 
  
  Now, for $\u\in\{0\}\times\Z^{s}$, let  $\W_{\u}(\omega)=\{\z\in\Z^{1+s}:\z\updownarrow \u\mbox{ in }\omega\}.$  
  We proceed to construct an auxiliary site percolation process of \textit{good} vertices on $\{0\}\times\Z^{s}$.  Given $\epsilon>0$ and a configuration $\omega$ we declare 
  each vertex $\u\in \{0\}\times\Z^{s}$ as \textit{good} if 
  \begin{enumerate}
    \item $\W_{\u}(\omega)>\frac{1+p}{1-p}\epsilon$
    \item there exists at least one open edge with exactly one endvertex in $\W_{\u}$ 
    in each of the $s$ possible directions of increasing coordinate-value in configuration $\omega$.
   \end{enumerate}
   
  Consider the sequence of independent events $\{A_{\u}(\omega)\}_{\u\in \{0\}\times\Z^{s}}$, where $$A_{\u}(\o)=\{\o\in\Omega:\u \mbox{ is \textit{good} in the configuration }\omega\},\,\,\,\,\,\u\in \{0\}\times\Z^{s}.$$ 
We aim to obtain an estimate for the probability of $A_{\u}$. If we denote the event in condition 2. by $F_{\u}$, then for any $\u\in \{0\}\times\Z^{s}$, we have
 \begin{equation}\label{prob_good}
   \mathbb{P}_{p,q}(\u\mbox{ is good})=\mathbb{P}_{p,q}\left(\W_{\u}>\frac{1+p}{1-p}\epsilon;\,F_{\u}\right).
   \end{equation}
 Now, for $\epsilon>0$ sufficiently small and $p$ sufficiently close to 1,
 
 \begin{equation}\label{cluster}
   \mathbb{P}_{p,q}\left(\W_{\u}>\frac{1+p}{1-p}\epsilon\right)=p^{\frac{1+p}{1-p}\epsilon}\geq 
   1-3\epsilon.
 \end{equation}


Define $$F_{\u,j}(\o)=\{\mbox{ there exists at least one open edge with exactly one endvertex in $\W_{\u}$} $$
$$\mbox{ in direction $j$ in configuration $\o$}\},\,\,\,j=1,\dots,s.$$ Note that  $\{F_{\u,j}\}_{j}$, $j=1,\dots,s$, is a collection of independent events with equal probabilities. Then
 $$\mathbb{P}_{p,q}\left(F_{\u}|\W_{\u}>\frac{1+p}{1-p}\epsilon\right)=\mathbb{P}_{p,q}\left(\bigcap_{j=1}^s F_{\u,j}\Big|\W_{\u}>\frac{1+p}{1-p}\epsilon\right)$$
   $$\geq \left[\mathbb{P}_{p,q}\left(F_{\u,1}\Big|\W_{\u}>\frac{1+p}{1-p}\epsilon\right)\right]^s\geq\left[1-(1-q)^{\frac{1+p}{1-p}\epsilon}\right]^s.$$
 If for some $\alpha>0$ we choose $q>\alpha\frac{1-p}{1+p}$,  then
  \begin{equation}\label{condicional}
   \mathbb{P}_{p,q}\left(F_{\u}\Big|\W_{\u}>\frac{1+p}{1-p}\epsilon\right)\geq\left[1-\left(1-\alpha\frac{1-p}{1+p}\right)^{\frac{1+p}{1-p}\epsilon}\right]^s\geq(1-e^{-\alpha\epsilon})^s.
 \end{equation}
  Plugging Equations (\ref{cluster}) and (\ref{condicional}) into Equation (\ref{prob_good}) we obtain
\begin{equation*}\label{prob_origingood}
   \bar{p}:=\mathbb{P}_{p,q}(\u\mbox{ is good})\geq 
   (1-3\epsilon)(1-e^{-\alpha\epsilon})^s.
 \end{equation*}
 Taking $\alpha=\alpha(\epsilon)$ sufficiently large, this can be made strictly 
 larger than the critical threshold of site percolation on $\Z^s$. Since the anisotropic bond process with parameters $(p,q)$ stochastically dominates the isotropic site process with parameter $\bar{p}$, we obtain
 $$\mathbb{P}_{p,q}(\mbox{ bond percolation occurs in }\Z^{1+s})\geq \mathbb{P}_{\bar{p}}(\mbox{ site percolation occurs in }\Z^{s} )>0.$$ This completes the proof.

\section*{Acknowledgements} R\'emy Sanchis was partially supported by Conselho Nacional de Desenvolvimento Cient\'ifico e Tecnol\'ogico (CNPq) and by  Fundac\~ao de Amparo \`a Pesquisa do Estado de Minas
Gerais (FAPEMIG), grant PPM 00600/16. Roger W. C. Silva was partially supported by FAPEMIG, grant APQ-02743-14.


\begin{thebibliography}{}
\bibitem{BK} Berg J. van den, Kesten H., Inequalities with applications to percolation and reliability. {\em  Journal of Applied Probability.}; {\bf{22}}: 556-569, 1985.

\bibitem{BMT} Blanc R., Mitescu C.D, and  Th\'{e}venot G.,  Percolation anisotrope : conductivit\'{e} dâun r\'{e}seau carr\'{e} de liens al\'{e}atoires. {\em Journal de Physique}; {\bf {41}}: 387-391, 1980.

\bibitem{CLS} Couto R.G., de Lima, B.N.B. and Sanchis, R.,  Anisotropic percolation on slabs.
{\em  Markov Processes and Related Fields.}; {\bf{20}}: 145-154, 2014.

\bibitem{DB} Deng Y. and  Bl\"{o}te H. W. J., Anisotropic limit of the bond-percolation model and conformal invariance in curved geometries. {\em Physical Review E}; {\bf {69}}: 066129,  2004.

\bibitem{FMMPV} Fontes L.R., Marchetti, D.H.U., Merola, I., Presutti, E. and Vares, M.E.,  Phase transitions in layered systems.
{\em  Journal of Statistical Physics.}; {\bf{157}}: 407-421, 2014.

\bibitem{FMMPV2} Fontes L.R., Marchetti, D.H.U., Merola, I., Presutti, E. and Vares, M.E.,  Layered systems at the mean field critical temperature.
{\em  Journal of Statistical Physics.}; {\bf{161}}: 91-122, 2015.

\bibitem{FS} Friedman S. P. and Seaton N. A., Percolation thresholds and conductivities of a uniaxial anisotropic simple-cubic lattice. {\em Transport in Porous Media};  {\bf{30}}: 241-250, 1998.

\bibitem{GCGR} Guyon E., Clerc J.P., Giraud G. and Roussenq J.,   A network simulation of anisotropic percolation.
{\em  Journal de Physique.}; {\bf{42}}: 1553-1557, 1981.

\bibitem{GRI} Grimmett G., Percolation, 2nd ed., Springer-Verlag, Berlin, 1999. 

\bibitem{GN} Grimmett G. and Newman C.M.,  Percolation in $\infty+1$ dimensions. {\em  Disorder in Physical Systems} (ed. G.R. Grimmett and D.J.A. Welsh); Clarendon Press, Oxford: 219-240, 1990.

\bibitem{K} Kesten H., Percolation Theory for Mathematicians, Birkh{\"a}user, Boston, 1982. 

\bibitem{M}  McGurn A. R., Bond percolation on the square lattice: anisotropic systems. {\em Physlca A}; {\bf{119}}: 295-306, 1983.

\bibitem{MPS} Mazel A., Procacci, A. and Scoppola, B.,  Gas phase of asymmetric neighbor Ising model.
{\em  Journal of Statistical Physics.}; {\bf{106}}: 1241-1248, 2002.

\bibitem{RC}  Redner S. and Coniglio A.,  On the crossover exponent for anisotropic bond percolation. {\em Physics Letters};  {\bf{79}}: 111-112, 1980.

\bibitem{RS} Redner S. and Stanley H.E.,   Anisotropic bond percolation.
{\em  J. Phys. A: Math. Gen.}; {\bf{12}}: 1267-1283, 1979.

\bibitem{S} Stauffer D., Scaling properties of percolation clusters.
{\em  Disordered Systems and Localization}; Springer, Berlin, 9-25, 1981.











\end{thebibliography}
\end{document}